\documentclass{scrartcl}
\usepackage[utf8]{inputenc}
\usepackage{lmodern} 
\usepackage{amsmath}
\usepackage{amssymb}
\usepackage{amsthm}
\usepackage{thmtools}
\usepackage{mathtools}
\usepackage{hyperref} 
\usepackage{enumitem} 
\usepackage{graphicx} 
\usepackage{dsfont} 
\usepackage{tikz} 
\usepackage[framemethod=TikZ]{mdframed}
\usepackage{ifthen}
\usepackage{supertabular} 

\renewcommand{\subset}{\subseteq}
\newcommand{\lcb}{\left\lbrace} 
\newcommand{\rcb}{\right\rbrace} 
\newcommand{\cb}[1]{\lcb #1 \rcb} 
\newcommand{\lb}{\left(} 
\newcommand{\rb}{\right)} 
\newcommand{\br}[1]{\lb #1 \rb} 
\newcommand{\brOf}[1]{\!\br{#1}} 
\newcommand{\abs}[1]{\left| #1 \right|} 
\newcommand*{\E}{\mathbb{E}} 
\newcommand*{\V}{\mathbb{V}} 
\let\Pr\relax
\newcommand*{\Pr}{\mathbb{P}} 
\newcommand{\sizedMid}[2]{#1 \, \kern-\nulldelimiterspace\mathopen{}\left| \vphantom{#1}\,#2\right.\mathclose{}\kern-\nulldelimiterspace}

\newcommand{\Vof}[1]{\V[#1]}

\newcommand{\PrOf}[1]{\Pr\mathopen{}\lb #1 \rb\mathclose{}}
\newcommand{\Prof}[1]{\Pr(#1)}

\DeclarePairedDelimiterX\Set[1]{\lbrace}{\rbrace}%
{  #1 }
\newcommand{\Ex}{\E\expectarg}
\DeclarePairedDelimiterX{\expectarg}[1]{[}{]}{%
	\ifnum\currentgrouptype=16 \else\begingroup\fi
	\activatebar#1
	\ifnum\currentgrouptype=16 \else\endgroup\fi
}
\newcommand{\innermid}{\nonscript\;\delimsize\vert\nonscript\;}
\newcommand{\activatebar}{%
	\begingroup\lccode`\~=`\|
	\lowercase{\endgroup\let~}\innermid 
	\mathcode`|=\string"8000
}
\newcommand*{\mc}[1]{\mathcal{#1}}
\newcommand*{\mb}[1]{\mathbb{#1}}

\newcommand*{\ms}[1]{\mathsf{#1}}
\newcommand*{\mo}[1]{\mathbf{#1}}

\newcommand{\N}{\mathbb{N}}

\newcommand{\R}{\mathbb{R}}

\newcommand{\pr}{^\prime}

\def\integral from #1to #2of #3by #4;{\int_{#1}^{#2} \! #3 \mathrm{d}#4} %
\def\integralMeasure in #1of #2by #3of #4;{\int_{#1} \! #2{#4} #3{\mathrm{d}#4}} %
\def\mapping #1from #2to #3;{#1 \colon #2 \rightarrow #3}
\def\mappingDef #1from #2to #3maps #4to #5;{#1 \colon #2 \rightarrow #3,\ #4 \mapsto #5}
\def\seq #1by #2;{\br{#1}_{#2\in\N}}
\def\seqInText #1by #2;{(#1)_{#2\in\N}}
\newcommand{\lebesgue}{\mathcal{L}}
\newcommand{\lebesguePow}[1]{\lebesgue^{#1}}

\newcommand{\sgn}{\mathsf{sgn}} 
\newcommand{\dl}{\mathrm{d}}
\def\converges for #1to #2;{\xrightarrow{#1} #2}
\def\convergesAlmostSurely for #1to #2;{\xrightarrow{#1}_{\mathsf{fs}} #2}
\def\convergesInProbability for #1to #2;{\xrightarrow{#1}_{\mathsf{p}} #2}
\def\convergesInL #1for #2to #3;{\xrightarrow{#2}_{\lebesguePow{#1}} #3}

\newcommand{\ind}{\mathds{1}}
\newcommand{\normof}[1]{\Vert #1 \Vert}
\newcommand{\normOf}[1]{\left\Vert #1 \right\Vert}
\newcommand{\equationFullstop}{\, .}
\newcommand{\eqfs}{\equationFullstop}
\newcommand{\equationComma}{\, ,}
\newcommand{\eqcm}{\equationComma}

\newcommand{\euler}{\mathrm{e}}
\DeclareMathOperator*{\argmin}{arg\,min}

\def\postBoxSkip{1.0ex}
\def\postBoxSkipCmd{\vskip\postBoxSkip}
\def\preBoxSkip{1.0ex}
\def\preBoxSkipCmd{\vskip\preBoxSkip}
\declaretheoremstyle[
	bodyfont=\normalfont,
	postfoothook={\postBoxSkipCmd},
	preheadhook={\preBoxSkipCmd},
	mdframed={
		backgroundcolor = black!2,
		startcode={\def\environmentEnumerateLabel{(\roman*)}},
}]{ruledBoxStyle}
\declaretheoremstyle[
	bodyfont=\normalfont,
	postfoothook={\postBoxSkipCmd},
	preheadhook={\preBoxSkipCmd},
	mdframed={
		backgroundcolor=white,
}]{ruledBoxStyleWhite}
\declaretheoremstyle[
	bodyfont=\normalfont,
	postfoothook={\postBoxSkipCmd},
	preheadhook={\preBoxSkipCmd},
	mdframed={
		backgroundcolor=black!2,
		linecolor = black!2,
		tikzsetting = {
			draw = black,
			line width = 2pt,%
			dashed,%
			dash pattern = on 10pt off 3pt
		},
}]{dashedBoxStyle}
\declaretheoremstyle[
	bodyfont=\normalfont,
	headformat={\NAME \NOTE},
	postfoothook={\postBoxSkipCmd},
	preheadhook={\preBoxSkipCmd},
	mdframed={
		linecolor = white,
		startcode={\def\environmentEnumerateLabel{(\roman*)}},
		tikzsetting = {
			draw = black,
			line width = 1pt,%
			loosely dotted,
		},
	}
]{dashedStyle2}
\declaretheoremstyle[
	bodyfont=\normalfont,
	postfoothook={\postBoxSkipCmd},
	preheadhook={\preBoxSkipCmd},
	mdframed={
		linecolor = white,
		startcode={\def\environmentEnumerateLabel{(\roman*)}},
		tikzsetting = {
			draw = black,
			line width = 1pt,%
			loosely dotted,
		},
	}
]{dashedStyle}
\declaretheoremstyle[
	bodyfont=\normalfont,
	postfoothook={\postBoxSkipCmd},
	preheadhook={\preBoxSkipCmd},
	mdframed={
		linecolor = black,
		innerlinewidth=1pt,outerlinewidth=1pt,
		middlelinewidth=1pt,
		linecolor=black,middlelinecolor=white,
		startcode={\def\environmentEnumerateLabel{(\roman*)}},
	}
]{doubleStyle}
\declaretheoremstyle[
	bodyfont=\normalfont,
	postfoothook={\postBoxSkipCmd},
	preheadhook={\preBoxSkipCmd},
	mdframed={
		backgroundcolor = black!4,
		linecolor = black!4,
		startcode={\def\environmentEnumerateLabel{(\alph*)}},
}]{boxStyle}
\declaretheoremstyle[
	headfont=\normalfont\itshape, 
	notefont=\normalfont\itshape, 
	notebraces={}{},
	bodyfont=\normalfont,
	qed=\qedsymbol,
	numbered=no,
	headindent=0pt,
	postheadspace=1ex,
	name={Proof},
	postheadhook={\def\environmentEnumerateLabel{(\roman*)}},
	mdframed={
		hidealllines = true,
		innerrightmargin = 0pt,
		innerleftmargin = 0pt,
		innertopmargin = 0pt,
		innerbottommargin = 0pt,
		leftmargin = 0pt,
		rightmargin = 0pt,
	}
]{proofStyle}
\declaretheoremstyle[
	bodyfont=\normalfont,
	postfoothook={\postBoxSkipCmd},
	preheadhook={\preBoxSkipCmd},
	mdframed={
		backgroundcolor = white,
		linecolor = black,
		startcode={\def\environmentEnumerateLabel{(\alph*)}},
		leftline = false,
		rightline = false,
}]{tobBottomStyle}
\declaretheoremstyle[
bodyfont=\normalfont,
]{standardStyle}
\declaretheorem[style=ruledBoxStyle,name=Definition]{definition}
\declaretheorem[style=ruledBoxStyle,name=Lemma,numberlike=definition]{lemma}
\declaretheorem[style=ruledBoxStyle,name=Proposition,numberlike=definition]{proposition}

\declaretheorem[style=ruledBoxStyle,name=Theorem,numberlike=definition]{theorem}
\declaretheorem[style=ruledBoxStyle,name=Corollary,numberlike=definition]{corollary}
\declaretheorem[style=ruledBoxStyle,name=Theorem,numbered=no]{theorem*}
\declaretheorem[style=boxStyle,name=Remark,numberlike=definition]{remark}

\makeatletter
\let\proof\@undefined
\let\endproof\@undefined
\makeatother
\declaretheorem[style=proofStyle]{proof}
\def\theoremContentInNewLine{\text{}}
\def\environmentEnumerateLabel{(\roman*)}
\newcounter{subExample}%
\renewcommand{\thmcontinues}[1]{Teil \arabic{subExample}} 
\newcommand{\myQf}{\mc Q}
\newcommand{\myg}{g}
\newcommand{\rvs}{,}
\begin{document}
\title{Arbitrary Rates of Convergence for Projected and Extrinsic Means}
\author{Christof Schötz\thanks{Institute of Applied Mathematics,
Heidelberg University,
Im Neuenheimer Feld 205,
69120 Heidelberg, Germany,
schoetz@math.uni-heidelberg.de}}
\maketitle
\begin{abstract} 
We study central limit theorems for the projected sample mean of independent and identically distributed observations on subsets $\mathcal Q \subset \mathbb R^2$ of the Euclidean plane. 

It is well-known that two conditions suffice to obtain a parametric rate of convergence for the projected sample mean:  $\mathcal Q$ is a $\mathcal C^2$-manifold, and the expectation of the underlying distribution calculated in $\mathbb R^2$ is bounded away from the medial axis, the set of point that do not have a unique projection to $\mathcal Q$.

We show that breaking one of these conditions can lead to any other rate: For a virtually arbitrary prescribed rate, we construct $\mathcal Q$ such that all distributions with expectation at a preassigned point attain this rate.
\end{abstract}
\section{Introduction}
Let $Z$ be a random variable with values in $\R^2$ and finite second moment. Let $\mc Q \subset \R^2$ be a subset of the Euclidean plane. 
Assume $m = \argmin_{p\in\mc Q} \normof{\Ex{Z}-p}$ exists and is unique. We call $m$ the \textit{projected (population) mean} of $Z$ in $\mc Q$. 
Let $Z_1, \dots, Z_n$ be independent and identically distributed copies of $Z$. We estimate $m$ by a \textit{projected sample mean} $m_n \in \argmin_{p\in\mc Q} \normof{\bar Z_n-p}$, $\bar Z_n := \frac1n\sum_{i=1}^n Z_i$. If $Z$ takes values only in $\mc Q$, then $m$ and $m_n$ are called \textit{extrinsic (population) mean} and \textit{extrinsic sample mean}, respectively \cite{hendriks98, bhattacharya03}. In \cite{hendriks98}, the extrinsic mean is called mean location.

For a given rate sequence $(a_n)_{n\in\N} \subset(0, \infty)$, $a_n\to0$ our goal is to find a set $\mc Q$ such that for a large class of distributions of $Z$ a central limit theorem of the form $a_n^{-1} (m_n - m) \xrightarrow{n\to\infty} \nu$ holds for some non-degenerate distribution $\nu$. Then $m_n$ converges to $m$ in probability at rate $a_n$.

Asymptotics of extrinsic sample means in cases with parametric rate of convergence, i.e., $a_n = n^{-\frac12}$ , are well-studied \cite{hendriks98, patrangenaru98, bhattacharya03, bhattacharya05}. This line of work is mostly concerned with finite dimensional manifolds, but results for infinite dimensional Hilbert manifolds are available \cite{ellingson13}. Slower rates for \textit{intrinsic} sample means, i.e., minimizers of $p\mapsto\sum_{i=1}^n d_{\mc Q}(Z_i, p)^2$ with the intrinsic metric $d_{\mc Q}$, have been observed on the circle \cite{hotz15} and more general manifolds \cite{huckemann18}. In some cases intrinsic and extrinsic means coincide \cite[Theorem 3.3]{bhattacharya03}. But this is not true in general.
 
The occurrence of a rate of convergence slower than the parametric one is called \textit{smeariness}. If, in contrast, the sample mean is equal to its population counterpart with high probability, the behavior is called \textit{stickiness}, which is observed for intrinsic means in certain negatively curved spaces \cite{hotz13, huckemann15}. 

\subsection{Medial Axis and Reach}

Our analysis is strongly connected to the \textit{medial axis} $\mc M_{\mc Q}$ of the set $\mc Q$, which is the set of all points that have more than one closest point in $\mc Q$. Formally,
\begin{equation*}
\mc M_{\mc Q} = \cb{z \in \R^2 \ \bigg\vert\  \exists p_1, p_2 \in \mc Q, p_1 \neq p_2\colon \normof{p_1 - z} = \normof{p_2 - z} = \inf_{p\in\mc Q} \normof{p - z}}\eqfs
\end{equation*}
The medial axis has been analyzed from a purely geometric perspective \cite{birbrair17}. The \textit{reach} \cite{federer59} $\tau_{\mc Q} := \inf_{m \in \mc M_{\mc Q}, p\in \mc Q} \normof{m-p}$ of a set $\mc Q \subset \R^2$ is the largest nonnegative real value such that any point in $\R^2$ with distance to $\mc Q$ less than $\tau_{\mc Q}$ has a unique closest point in $\mc Q$.

By the definition of medial axis $\mc M_{\mc Q}$ as it it is used here, it need not be a closed set, as the example $\mc Q = \cb{y = x^2} \subset \R^2$, $\mc M_{\mc Q} = (1/2, \infty) \times \cb{0}$ shows. Note that this contrasts some mentions of the term in the literature, e.g., in the context of \cite[Theorem 3.2]{bhattacharya03}. See \cite[Theorem A.5]{huckemann10} for a sufficient condition for a closed medial axis.

If $\mc Q$ is a $\mc C^2$-manifold, the projection map $z \mapsto \Pi_{\mc Q}(z) = \argmin_{p \in \mc Q} \normof{z - p}$ is continuously differentiable on $\R^2 \setminus \overline{\mc M_{\mc Q}}$ with $\normof{\nabla \Pi_{\mc Q}(z)} > 0$ \cite{abatzoglou78}. If additionally the reach $\tau_{\mc Q}$ is greater than distance of $\Ex{Z}$ to $\mc Q$, then the projected sample mean attains a parametric rate of convergence \cite{hendriks98, bhattacharya05}: The delta-method yields $\sqrt{n}\br{m_n - m} = \sqrt{n}\br{\Pi_{\mc Q}(\bar Z_n) - \Pi_{\mc Q}(\Ex{Z})} \xrightarrow{d} \mc N(0 , \tilde \Sigma)$ where $\tilde \Sigma = \nabla \Pi_{\mc Q}(\Ex{Z})\pr \cdot \mb{COV}(Z) \cdot \nabla\Pi_{\mc Q}(\Ex{Z})$. As convergence is a local phenomenon, we can replace the condition on the reach by the requirement that $\Ex{Z}$ is bounded away from the medial axis $\mc M_{\mc Q}$.

We construct sets $\mc Q$ with faster and slower rates of convergence than $1/\sqrt{n}$.
In our examples, the sets $\mc Q$ for slow rates are $\mc C^2$-smooth but $\Ex{Z}$ is too close to the medial axis, i.e., $\Ex{Z} \in \overline{\mc M_{\mc Q}}\setminus \mc M_{\mc Q}$. Sets $\mc Q$ with fast rates have reach $\tau_{\mc Q} > \inf_{p\in\mc Q} \normof{\Ex{Z} - p}$ but are only $\mc C^1$- but not $\mc C^2$-manifolds.

\subsection{Our Construction}

\begin{figure}
	\includegraphics[width=\textwidth]{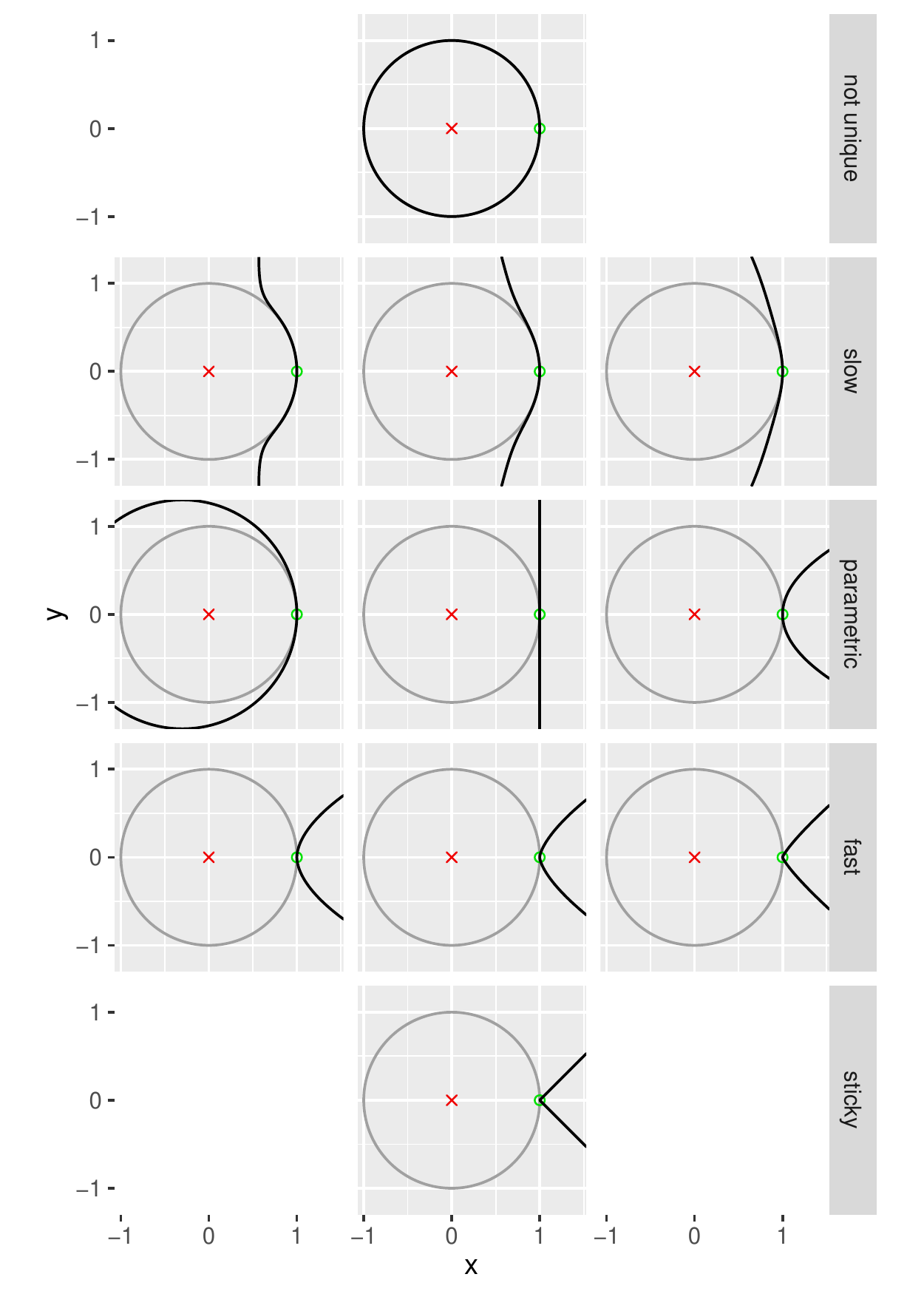}
	\caption{The images show the transition of the set $\mc Q$ (black) from non-unique projections, to slow, parametric, and fast rates, and sticky behavior of $m_n$. For reference, a circle (gray) with radius 1 around the origin is drawn. The expectation of $Z$ and its projection to $\mc Q$ are marked in red and green, respectively.}%
	\label{fig:quali}
\end{figure}

For a continuous function $f$ with $f(0) = 0$, we construct $\mc Q = \mc Q_f$ such that the projection of a point $(x\rvs y)\pr\in\R^2$ to $\myQf$ is roughly $(1 \rvs f(y))\pr$ for $|x|, |y|$ small enough. Assuming $\Ex{Z} = 0 \in\R^2$, the arithmetic mean $\bar Z_n = (\bar X_n\rvs \bar Y_n)\pr$ concentrates at 0 with rate $1/\sqrt{n}$. Thus, $m_n = \Pi_{\mc Q}(\bar Z_n) \approx (1\rvs f(\bar Y_n))\pr$ concentrates at $(1 \rvs 0)\pr$ with a rate depending on $f$.
For a wide range of rates $(a_n)_{n \in \N} \subset(0, \infty)$, $a_n\to0$, we can find a function $f$ with corresponding set $\myQf$ such that $a_n^{-1} (m_n - m) \xrightarrow{n\to\infty} \nu$ in distribution for some non-degenerate distribution $\nu$. 
As an example, $f(y) = \abs{y}^\gamma$, $\gamma > 0$ yields $a_n = n^{-\frac\gamma2}$, see \autoref{coro}. 
Examples of the constructed sets for qualitatively different rates can be found in \autoref{fig:quali}.

\subsection{Outline}

In \autoref{sec:res}, we present our main theoretical results. We state the requirements on the function $f$ and describe how the set $\myQf$ is constructed from $f$. \autoref{prop} states a result on deterministic projection to $\myQf$, while \autoref{thm} describes how the projected sample mean converges to the projected population mean. The goal of \autoref{sec:illu} is to illustrate the general statement of \autoref{thm}. We first derive \autoref{coro}, which gives explicit functions $f$  and sets $\myQf$ for certain prescribed rates $a_n$. In particular, we give examples where projected means attain polynomial, logarithmic, or exponential rates of convergence. Then the results are discussed and visualized. All proofs are given in \autoref{sec:proof}.
\section{Results}\label{sec:res}
The possible choices of the function $f$, which determines the set $\mc Q = \mc Q_f$ and, thus, the rate of convergence, are not restricted very much.
\begin{itemize}
\item[\textbf{(A0)}:]
Let $b > 0$. Let $f \in \mc C^0([0,b])$ be strictly increasing with $f(0) = 0$. 
\end{itemize}
Under the assumption \textbf{(A0)}, we construct the set $\mc Q$ as follows. Set $B := f(b)$.
We denote the inverse function of $f\colon [0,b] \to [0, B]$ by $\myg\colon [0,B] \to [0, b]$, i.e., $\myg(x) := f^{-1}(x)$. For $t\in [0,B]$, define $r(t) = r_f(t) := 1 + \int_0^t \myg(x) \dl x$. 
Finally, define 
\begin{align}\label{eq:q}
\begin{split}
	q(t) = q_f(t) &:= r(\abs{t}) \begin{pmatrix} \cos(t) \\ \sin(t) \end{pmatrix} \text{ for } t \in [-B, B]\eqcm\\
	\mc Q = \mc Q_f &:= \cb{ q(t) \colon t \in [-B, B]}\eqfs
\end{split}
\end{align}
Our main results are based on the observation that the projection of a point $(x\rvs y)\pr$ to $\mc Q$ for $x, y$ small enough is essentially $(1\rvs f(y))\pr$.

We denote the projection of $z \in \R^2$ to $\mc Q$ as $\Pi_{\mc Q}(z)$, i.e., $\Pi_{\mc Q}(z) = \argmin_{p \in \mc Q} \normof{z - p}$. If the argmin is not unique, we assume that one element of the argmin--set is chosen by a fixed arbitrary mechanism, e.g., smallest lexicographic order. The argmin--set cannot be empty as $\mc Q$ is compact by construction.
\begin{lemma}\label{lmm}
Assume \textbf{(A0)} with $\mc Q$ from \eqref{eq:q}.
Let $y \in \R$ with $y \to 0$, $x = \mo O(y)$, and $t_y \in [-B, B]$ such that $\Pi_{\mc Q}((x \rvs y)\pr) = q(t_y)$.
Then
\begin{equation*}
	\myg(t_y) = y + \mo o(y)
	\eqfs
\end{equation*}
\end{lemma}
\begin{remark}[Simpler construction]
As can be seen from the proof of \autoref{lmm}, a simpler construction in the case of $f(t) = \mo o(t)$ is replacing $q(t)$ by
\begin{equation*}
	\tilde q(t) := \begin{pmatrix} 1 + t\myg(t) \\ t \end{pmatrix}
	\eqfs
\end{equation*}
This yields the same results, but it does not work for $\myg(t) = \mo o (t)$.
\end{remark}
\autoref{lmm} describes the projection of a point close to the origin in an indirect way, i.e, after applying the function $g$. To have a direct statement, we need to make additional assumptions.
\begin{itemize}
\item[\textbf{(A1)}:]
Assume
\begin{align*}
\lim_{y\searrow0} \frac{f\brOf{y + c y (y + f(y))}}{f(y)} = 1
\end{align*}
for all $c \in \R$.
\item[\textbf{(A1)'}:]
Assume
\begin{align*}
\lim_{y\searrow0} \frac{f\brOf{y + c y f(y) (y + f(y))}}{f(y)} = 1
\end{align*}
for all $c \in \R$.
\end{itemize}
\begin{proposition}\label{prop}
Assume \textbf{(A0)} and \textbf{(A1)} with $\mc Q$ from \eqref{eq:q}. 
Let $y \in \R$ with $y \to 0$, $x = \mo O(y)$, and $t_y \in [-B, B]$ such that $\Pi_{\mc Q}((x \rvs y)\pr) = q(t_y)$.
Then
\begin{equation*}
	t_y = f(y) + \mo o(f(y)) 
	\qquad\text{and}\qquad
	\Pi_{\mc Q}\brOf{ \begin{pmatrix} x \\ y \end{pmatrix}} =\begin{pmatrix} 1 \\ f(y) \end{pmatrix} + \mo o(f(y))
	\eqfs
\end{equation*}
Furthermore, if $x = 0$, we can replace the assumption \textbf{(A1)} by \textbf{(A1)'}.
\end{proposition}
\begin{remark}[On the assumptions \textbf{(A1)} and \textbf{(A1)'}]\label{rem1}
\theoremContentInNewLine
\begin{enumerate}[label=\environmentEnumerateLabel]
\item 
	We have
	\begin{equation}\label{eq:reltoone}
		\lim_{y\searrow0} \frac{f\brOf{y + \mo o(y)}}{f(y)} = 1
	\end{equation}
	for any function of the form $f(y) = ay^\gamma$, with $a,\gamma>0$.
	Furthermore, \eqref{eq:reltoone} implies \textbf{(A1)}, and \textbf{(A1)} implies \textbf{(A1)'}.
\item 
	It is unclear to the author, whether there is a function that fulfills \textbf{(A0)} but not \textbf{(A1')}.
\item	
	The function $f(y) = \exp(-1/y)$ fulfills \textbf{(A0)} and \textbf{(A1)'}, but does not fulfill \textbf{(A1)}. If we set $x = y$, we obtain \begin{equation*}
		t_y = f\brOf{\frac{y}{1-y}} + \mo o\brOf{f\brOf{\frac{y}{1-y}}} = \exp(1) f(y) + \mo o (f(y)) \neq  f(y) + \mo o (f(y))
		\eqfs
	\end{equation*}
	If we set $\tilde f(y) = \exp(-\exp(1/y))$, we even have $\tilde f(y) = \mo o(\tilde t_y)$.
	
	Note that $x \mapsto \exp(-1/x)\ind_{(0,\infty)}(x)$ is a classical example of a function that is infinitly often differentiable but not analytic: for every $k\in\N_0$ the $k$-th derivative at 0 is 0, $f^{(k)}(0) = 0$.
\item
	If $f \in \mc C^k$, i.e., $f$ is $k$-times continuously differentiable, $k \in \N$, and there is an $\ell \in \cb{1, \dots, k}$ such that $f^{(\ell)}(0) \neq 0$, we set $\ell_0 := \min \cb{\ell\in\cb{1, \dots, k} \colon f^{(\ell)}(0) \neq 0}$. Then, by Taylor's theorem, $f(z) = \frac{f^{(\ell_0)}(0)}{\ell_0!} z^{\ell_0} + \mo o (z^{\ell_0})$. Thus, \eqref{eq:reltoone}, \textbf{(A1)}, and \textbf{(A1)'} hold.
\end{enumerate}
\end{remark}
As taking the projected mean is projecting the Euclidean mean in $\R^2$ to $\mc Q$, \autoref{lmm} induces a central limit theorem for projected means.
\begin{theorem}\label{thm}
Assume \textbf{(A0)} with $\mc Q$ from \eqref{eq:q}.
Let $Z = (X \rvs  Y)\pr$ be a random variable in $\R^2$ with finite second moment, $\Ex{Z} = 0 \in \R^2$, and $\Vof{Y} = \sigma^2 > 0$. 
Let $Z_1, \dots, Z_n$ be independent copies of $Z$.
Then the projected population mean $m \in \mc Q$ exists, is unique, and
\begin{equation*}
	m = \Pi_{\mc Q} \brOf{\Ex{Z}} = \argmin_{p\in\mc Q} \Ex{\normof{Z - p}^2} = q(0) = \begin{pmatrix} 1 \\ 0 \end{pmatrix}
	\eqfs
\end{equation*}
Let $(m_{n,1}\rvs m_{n,2})\pr := m_n := \Pi_{\mc Q} \brOf{\bar Z_n}$, $\bar Z_n := \frac1n \sum_{i=1}^n Z_i$. Then $m_n$ is a projected sample mean. Let $t_n \in [-B,B]$ such that $m_n = q(t_n)$. Then, for $s\geq 0$,
\begin{align*}
	\lim_{n\to \infty} \PrOf{t_n \leq f\brOf{\frac{s}{\sqrt{n}}}} = 
	\lim_{n\to \infty} \PrOf{-t_n \leq f\brOf{\frac{s}{\sqrt{n}}}} &= \\
	\lim_{n\to \infty} \PrOf{m_{n,2} \leq f\brOf{\frac{s}{\sqrt{n}}}} =
	\lim_{n\to \infty} \PrOf{-m_{n,2} \leq f\brOf{\frac{s}{\sqrt{n}}}} &= 
	\Phi\brOf{\frac s\sigma}
	\eqcm
\end{align*}
where $\Phi$ denotes the distribution function of a standard normal random variable.
Moreover,
\begin{align*}
	\PrOf{ \abs{m_{n,1}-1} \geq f\brOf{\frac{s}{\sqrt{n}}}} 
	&\xrightarrow{n\to\infty}
	0\eqfs
\end{align*}
\end{theorem}
\begin{remark}[Arc length]\label{rem:arc}
	The curve $q(t)$ in \eqref{eq:q} is not necessarily parameterized by arc length. But $q \in \mc C^1((-B,B))$ and $\normof{\dot q(t)} = 1 + \mo o(1)$ for $\abs{t} \to 0$ as
	\begin{align*}
		\normof{\dot q(t)}^2 &= \br{\myg(t) \cos(t) - r(t)\sin(t)}^2 +  \br{\myg(t)\sin(t) + r(t)\cos(t)}^2
		\\&= 1 + (\myg(t)-t)^2 + \mo O(t\myg(t)+t^2)
		\eqfs
	\end{align*} 
	Thus, the results on $t_n$ in \autoref{thm} also hold if $t_n$ is replaced by an arc length parametrization.
\end{remark}
\begin{remark}[Why \autoref{thm} does not require \textbf{(A1)}]
In contrast to \autoref{prop}, we do not require \textbf{(A1)} or \textbf{(A1')} in \autoref{thm}. In particular, in the setting of \autoref{rem1} (c), $f(y) = \exp(-1/y)$, we have
\begin{equation*}
	\PrOf{t_n \leq f\brOf{\frac{s}{\sqrt{n}}}} \xrightarrow{n\to\infty} \Phi\brOf{\frac s\sigma}\eqcm
\end{equation*}
for $s \geq 0$ even though $t_n \neq  f(\bar Y_n) + \mo o (f(\bar Y_n))$. The reason is that the difference between $t_n$ and $f(\bar Y_n)$ is negligible in the scale that is used in \autoref{thm}. The right scale for a central limit theorem of $t_n$ is the one of $\bar Y_n$ (multiplied by $\sqrt{n}$), i.e., $\myg(t_n)$. The factor $\euler$ in $t_n \approx \euler f(\bar Y_n)$, see \autoref{rem1} (c), is non-negligible on the scale of $t_n$, but on the scale of $\bar Y_n$ it becomes 
\begin{equation*}
\frac{\myg(\euler f (\bar Y_n))}{\bar Y_n} = \frac{\log\brOf{\euler^{-1} \cdot \exp\brOf{\bar Y_n^{-1}}}^{-1}}{\bar Y_n} = 1/(1 - \bar Y_n) \xrightarrow{n\to\infty} 1
\end{equation*} almost surely, i.e., negligible.
\end{remark}
\begin{remark}[Non-uniqueness]
Non-unique closest points are not a problem in \autoref{thm} as $\Prof{\bar Z_n \in \mc M_{\mc Q}} \to 0$ by $\Vof{Y} = \sigma^2 >0$. See also \autoref{rem:reach}.
\end{remark}
\section{Illustration}\label{sec:illu}
\begin{figure}
	\includegraphics[width=\textwidth]{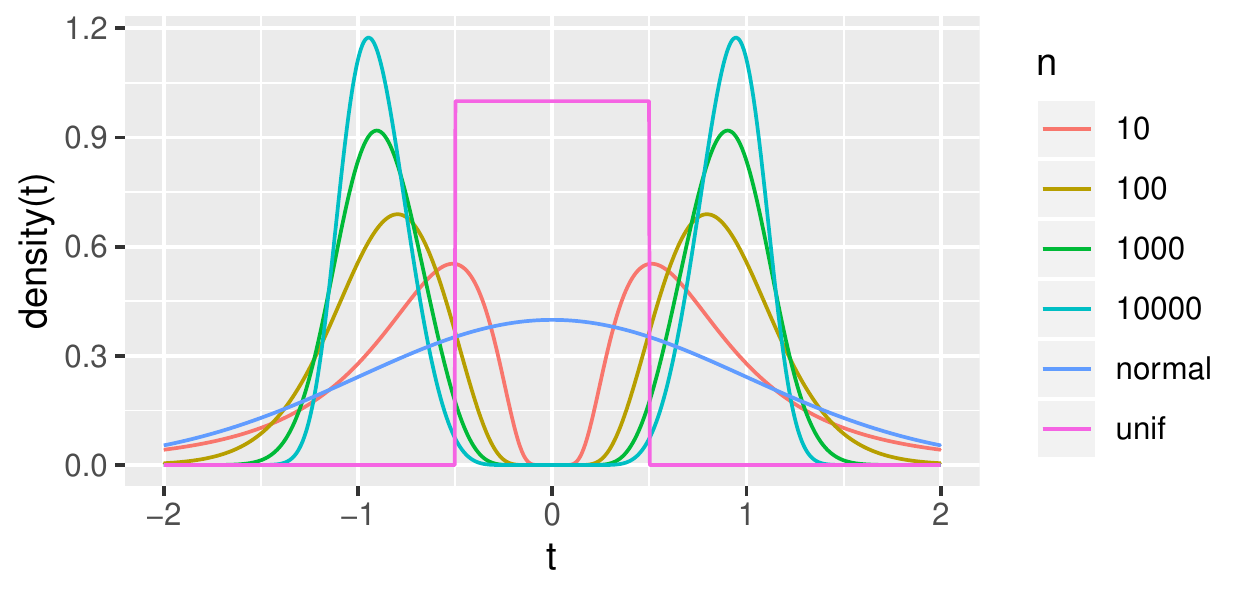}
	\caption{Plot of densities of $\frac12 \log(n) t_n$ for $f(y) = -\log(y)^{-1}$, $Z = (0\rvs Y)\pr$ and $Y \sim \mc N(0, 1)$, with standard normal and uniform densities for reference.}\label{fig:log_illu}
\end{figure}

\begin{figure}
	\includegraphics[width=\textwidth]{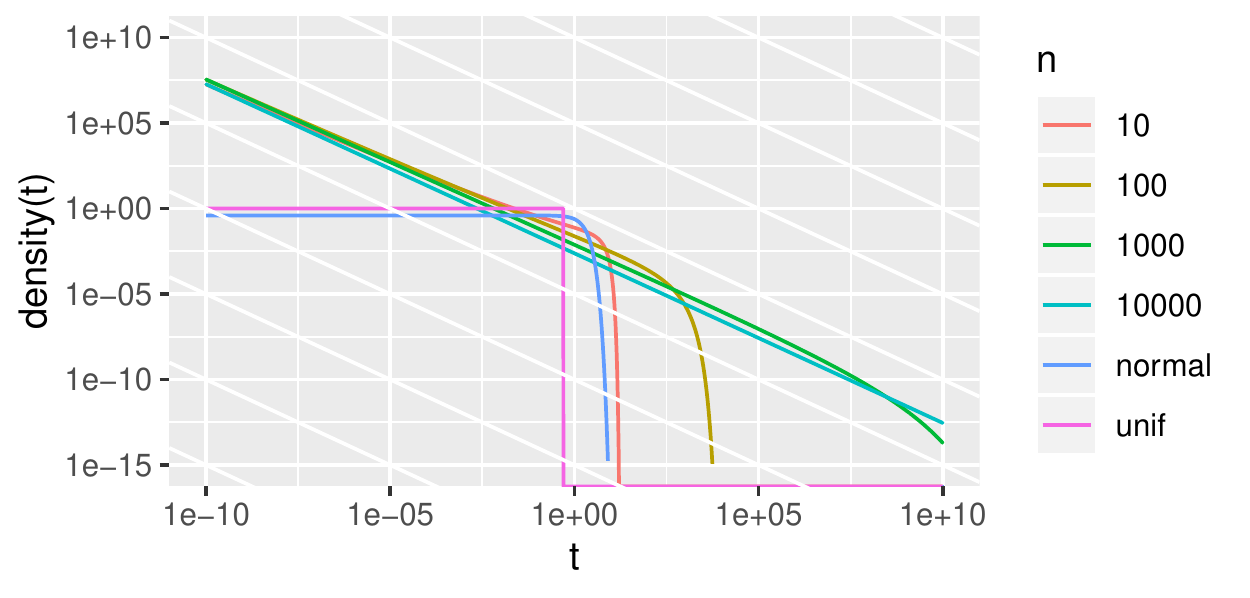}
	\caption{Log-log-plot of densities of $\exp(\sqrt{n}) t_n$ for $f(y) = \exp(-y^{-1})$, $Z = (0\rvs Y)\pr$ and $Y \sim \mc N(0, 1)$, with standard normal and uniform densities for reference.}\label{fig:exp_illu}
\end{figure}

To illustrate \autoref{thm}, we apply it to explicit functions $f$, which yield polynomial, logarithmic, and exponential rates of convergence for $m_n \to m$, respectively.
\begin{corollary}\label{coro}
Use the setting of \autoref{thm}.
	\begin{enumerate}[label=\environmentEnumerateLabel]
	\item 
	Let $f(y) = y^\gamma$ with $\gamma > 0$. Then $r(t) =  1+\frac{\gamma}{1 + \gamma} t^{\frac{1 + \gamma}{\gamma}}$ and
	\begin{equation*}
		n^{\frac{\gamma}{2}} t_n \to T
	\end{equation*}
	in distribution, where $\PrOf{T \leq s} = \Phi\brOf{\frac{\sgn(s)\abs{s}^\frac1\gamma}{\sigma}}$ for all $s \in \R$.
	\item  
	Let $f(y) = \br{-\log(y)}^{-\gamma}$ with $\gamma>0$. Then $r(t) =  1+\int_0^t\exp\brOf{-x^{-\frac1\gamma}} \dl x$ and 
	\begin{equation*}
		\br{\frac12\log(n)}^\gamma t_n \to T
	\end{equation*}
	in distribution,
	where $\Prof{T = 1} = \Prof{T = -1} = \frac12$.
	\item  
	Let $f(y) = \exp\brOf{-y^{-\gamma}}$  with $\gamma>0$. Then $r(t) =  1+\int_0^t\log\brOf{x^{-1}}^{-\frac1{\gamma}} \dl x$. For $c > 0$, define $U_{n,c} := \exp\brOf{(\sqrt{n}/c)^\gamma} t_n$ and $p_c := \Phi(\frac c\sigma)$.
	Then, for all $u\in(0,\infty)$, $\PrOf{U_{n,c} \geq u} \xrightarrow{n\to\infty} 1-p_c$, $\PrOf{U_{n,c} \leq -u} \xrightarrow{n\to\infty} 1-p_c$, and $\PrOf{\abs{U_{n,c}} \leq u} \xrightarrow{n\to\infty} 2p_c-1$.
	\end{enumerate}
	The results also hold when $t_n$ is replaced by $m_{n,2}$.
\end{corollary}
The results of \autoref{coro} are also true in arc length, see \autoref{rem:arc}.
\begin{remark}[On \autoref{coro}]
	\theoremContentInNewLine
	\begin{enumerate}[label=(\roman*)]
\item 
	For any polynomial scale $n^\gamma$, part (i) of \autoref{coro} gives an example of a central limit theorem with that scale.
\item  
	In part (ii) we obtain a central limit theorem with logarithmic scale and a Bernoulli-type limiting distribution that does not depend on $\sigma$. This seems quite remarkable and can be explained as follows:
	
	Scaling our observations $Z_i$ by $\sigma^{-1}$, is roughly like scaling $n$ by $\sigma^{2}$ as $\Vof{\sigma^{-1} \bar Y_n} = n^{-1} \approx \Vof{\bar Y_{[n\sigma^{2}]}}$, where $[n\sigma^{2}]$ denotes the closest integer to $n\sigma^{2}$. The scaling factor $\log(n)^\gamma$ is asymptotically equivalent to $\log(n\sigma^2)^\gamma$. Thus, constant factors like $\sigma$ cannot influence the asymptotic distribution on the scale $\log(n)^\gamma$.
	
	Densities of $t_n$ in the case of normally distributed observations are plotted in \autoref{fig:log_illu}.
\item  
	The statement of part (iii) of \autoref{coro}, can be summarized informally by
	\begin{equation*}
		\exp\brOf{(\sqrt{n}/c)^\gamma} t_n \to T_c
		\eqcm
	\end{equation*}
	where $\Prof{T_c = \infty} = \Prof{T_c = -\infty} = 1-p_c$ and $\Prof{T_c = 0} = 2p_c-1$.
	The limiting distribution has mass only at $0$ and $\pm\infty$. If the scale is changed such that the limit does not have a point mass at 0, all mass escapes to $\pm\infty$. If the scale is such that no mass escapes to $\pm\infty$, then in the limit all mass is at 0.
	
	Densities of $t_n$ in the case of normally distributed observations are plotted in \autoref{fig:exp_illu} on a log-log-scale. Only the positive axis of the symmetric densities is displayed.	
	The plot shows that the densities have non-negligible mass at all small orders of magnitude. Thus, choosing one specific order of magnitude by a specific scale makes all mass on larger orders of magnitude escape to infinity and all mass at smaller orders of magnitude go to 0.
%
	\end{enumerate}
\end{remark}

\begin{remark}[Extrinsic mean]
For the sets $\mc Q$ constructed in \autoref{coro}, there might not be a distribution with support in $\mc Q$ that has expectation 0. In particular, they might not directly yield examples of extrinsic means with the described asymptotic behavior. This is but a technical inconvenience. We can extend $\mc Q$ with an arbitrary set of points which have a distance to the origin that is bounded away from 1, and the result does not change. By doing so, we can also construct 2-dimensional manifolds with boundary which induce the same convergence results as the 1-dimensional structures in \autoref{coro}. 
\end{remark}
\begin{remark}[Application of \autoref{prop}]\label{rem2}
	For the functions $f$ in (i) and (ii), \textbf{(A1)} holds, see section \ref{ssec:proof:rem2}. Thus, \autoref{prop} implies 
	\begin{equation*}
		m_n = \Pi_{\mc Q}(\bar Z_n) \approx \begin{pmatrix}1\\\sgn(\bar Y_n)f(|\bar Y_n|)\end{pmatrix}\eqcm
	\end{equation*}
	meaning $\abs{m_{n,2} - \sgn(\bar Y_n)f(|\bar Y_n|)}/f(|\bar Y_n|) \to 0$ and $\abs{m_{n,1} - 1}/f(|\bar Y_n|) \to 0$ in probability.
	In (iii) only \textbf{(A1')} is true. Thus, the equation above is true for (iii) if $X = 0$ almost surely.
\end{remark}
\begin{remark}[Delta method]
In light of the delta method, note that, in the cases above, $f\pr(0)$ is 0 or $\infty$, except when $f$ is equal to the identity in (i). This is the only case of \autoref{coro} that yields the usual parametric rate. 
\end{remark}
\begin{figure}
	\includegraphics[width=\textwidth]{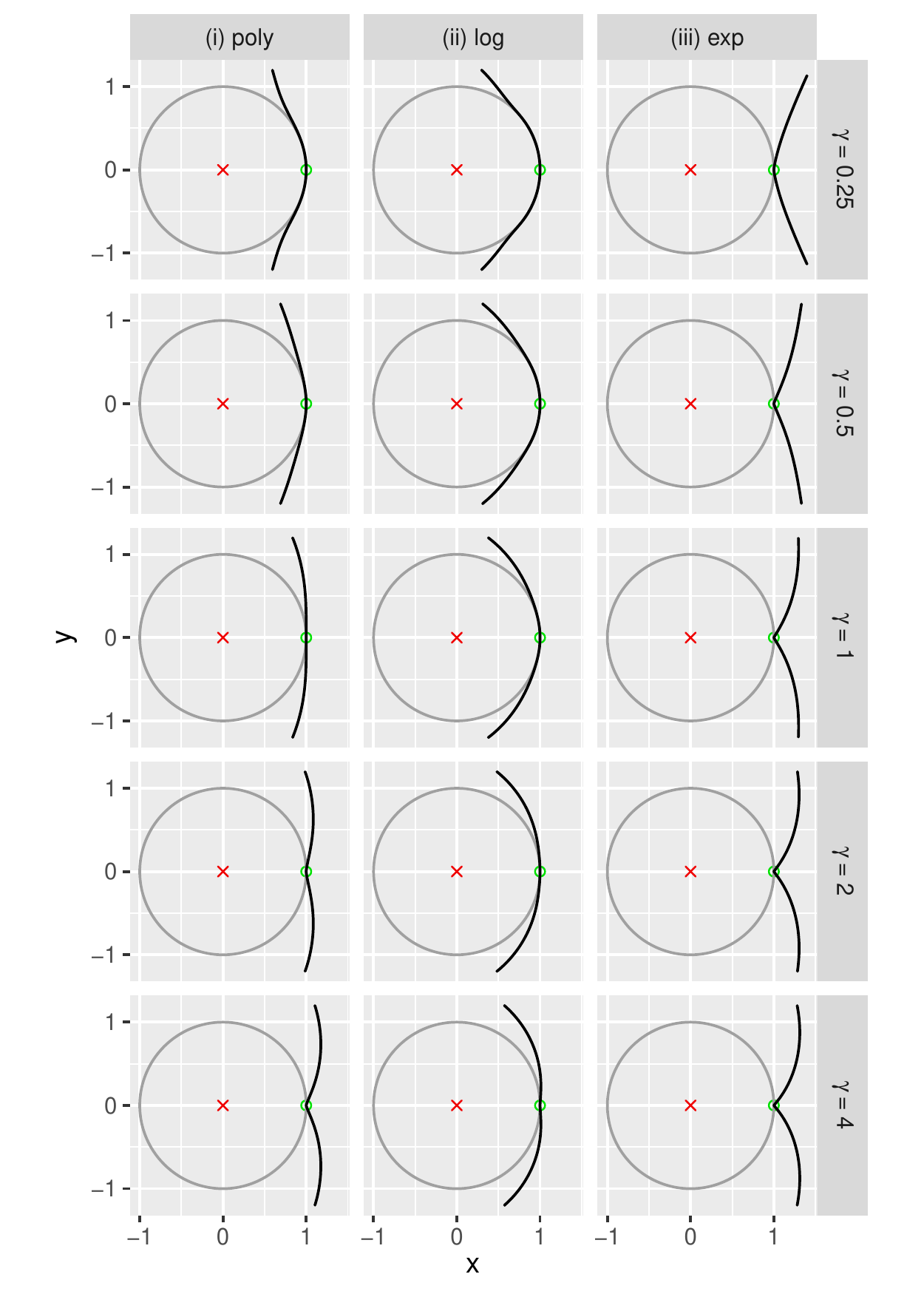}
	\caption{The images show the set $\mc Q$ (black) for different curves $q$, which are chosen as described in \autoref{coro}. For reference, a circle (gray) with radius 1 around the origin is drawn. The expectation of $Z$ and its projection to $\mc Q$ are marked in red and green, respectively.}\label{fig:full}
\end{figure}
Figure \ref{fig:full} illustrates the sets $\mc Q$ constructed according to the functions $f$ from \autoref{coro}.
The results on the convergence rate described in \autoref{thm} and \autoref{coro} depend only on the form of the curve close to the point $(1\rvs 0)\pr$. Even so the curve [\textit{(ii) log}, $\gamma=4$] looks like it is growing faster away from the circle than [\textit{(i) poly}, $\gamma=0.25$], the opposite is true when observing a neighborhood of $(1\rvs 0)\pr$ that is small enough.

\begin{figure}
	\centering
	\includegraphics[width=0.33\textwidth]{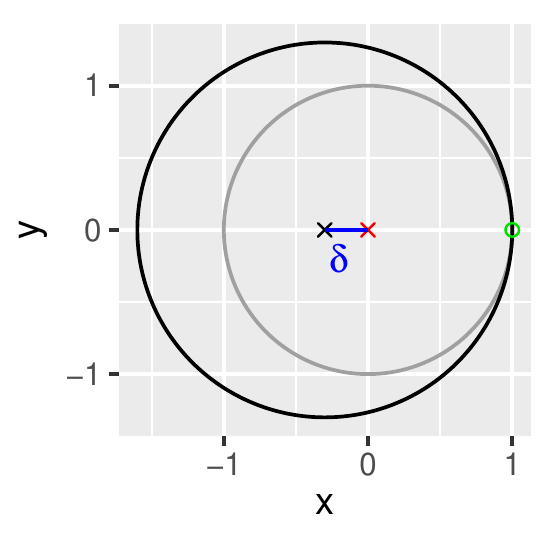}
	\caption{The black curve shows the set $\mc Q$ as described in \autoref{rem:circ} with $\delta = 0.3$. For reference, a circle (gray) with radius 1 around the origin is drawn. The expectation of $Z$ and its projection to $\mc Q$ are marked in red and green, respectively.}\label{fig:circ}
\end{figure}

There is a smooth transition of the set $\mc Q$ between slow and fast rates, see \autoref{fig:quali}. A circle with radius 1 centered at the origin can be seen as one extreme case, in the sense that an arbitrarily small change of a point at the origin can change its projection by a large amount. If $\mc Q$ almost looks like this circle, but increases its radius $r(t)$ slow enough, i.e, $r(t) \lesssim 1+t^2$, we still have large changes in the projection, but not arbitrarily large. For a larger circle with center $(-\delta\rvs 0)\pr$ and radius $1+\delta$ or a straight vertical line through $(1\rvs 0)\pr$ the changes of point and projection are proportional, i.e, $r(t) \approx 1 + t^2$. Changes in the point effect the projection only little if $q(t)$ grows to the right quickly when moving away from $(1\rvs 0)\pr$, i.e, $r(t) \gtrsim 1 + t^2$. For $\mc Q = \cb{(1\!+\!|y|\rvs  y)\pr\colon y\in \R}$ certain changes do not change the projection at all. In particular, $\Prof{m_n = m} \to 1$ (stickiness).
\begin{remark}[Larger circles]\label{rem:circ}
A circle with center at $(-\delta\rvs 0)\pr$, $\delta>0$, and radius $1+\delta$, see \autoref{fig:circ}, can be described by our construction with 
\begin{equation*}
	r(t) = \sqrt{\cos(t)^2 \delta^2 + 2\delta + 1} - \cos(t) \delta\eqcm
\end{equation*}
$t \in [-\pi, \pi]$.  Thus,
\begin{equation*}
	\myg(t) = \dot r(t) = \delta \sin(t) - \frac{\cos(t)\sin(t) \delta^2}{\sqrt{\cos(t)^2 \delta^2 + 2\delta + 1}} = \frac{\delta}{\delta+1} t + \mo O(t^2)\eqfs
\end{equation*}
Hence, the projection $\Pi_{\mc Q}(\bar Z_n)$ scales the $y$-direction only by a constant factor without affecting the rate of convergence. In particular we have a parametric rate of convergence. This can also be inferred by noting that $\mc Q$ is $\mc C^2$-smooth and has a reach larger than 1 as described in the introduction.
\end{remark}
\begin{remark}[Reach and Medial Axis]\label{rem:reach}
A set $\mc Q$ of our construction has reach at most 1 if $\myg(t) = \mo o(t)$ for $t \searrow 0$. This can be seen form \autoref{rem:circ}: If every circle with center at $(-\delta\rvs 0)\pr$ and radius $1+\delta$ for $\delta \in (0, \delta_0)$, $\delta_0 > 0$ intersects $\mc Q$ at more than one point the reach can be at most 1. Moreover, such a circle is constructed with $\myg_{\ms{circle}, \delta}(t)$ of order $t$, i.e., $\myg(t) = \mo o(\myg_{\ms{circle, \delta}}(t))$ and $r(t) = \mo o (r_{\ms{circle}, \delta}(t))$. Thus, $\cb{(-\delta\rvs 0)\pr \colon \delta \in (0, \delta_0)} \subset \mc M_{\mc Q}$ and $0 \in \partial \mc M_{\mc Q}$.
\end{remark}
\section{Proofs}\label{sec:proof}
\subsection{\autoref{lmm}}
Due to symmetry, we can restrict our analysis to $y \geq 0$ and $t \geq 0$ without loss of generality.
To find the projection point, we have to minimize the squared distance $\ell \in \mc C^1([-B,B])$, 
\begin{equation*}
	\ell(t) := \normOf{q(t) -  \begin{pmatrix}x \\ y \end{pmatrix}}^2 
	\eqfs
\end{equation*}
For its derivative, we have 
\begin{equation*}
	\frac12 \dot \ell(t) = r(t) \dot r(t) - x\br{\cos(t) \dot r(t) - \sin(t) r(t)} - y\br{\sin(t) \dot r(t) + \cos(t) r(t)}
	\eqfs
\end{equation*}
For $t\to 0$,
\begin{align*}
r(t) &= 1 + \mo O(t \myg(t))\eqcm\\
\dot r(t) &= \myg(t)\eqcm\\
\sin(t) &= t + \mo O(t^3)\eqcm\\
\cos(t) &= 1 + \mo O(t^2)\eqfs
\end{align*}
Thus, 
\begin{align*}
	\cos(t) \dot r(t) - \sin(t) r(t) &= \mo O(\myg(t) + t)\eqcm \\
	\sin(t) \dot r(t) + \cos(t) r(t) &= 1 + \mo O\brOf{t \myg(t) + t^2}\eqcm \\
	r(t) \dot r(t) &= \myg(t) + \mo O\brOf{t \myg(t)^2}\eqfs
\end{align*}
Denote by $t_y$ a global minimizer of $\ell(t)$. As $r(t)$ is strictly increasing for $t \geq 0$, we have $t_y \to 0$ as $x, y \to 0$.

Let $y \searrow 0$.
From $\dot \ell(t_y) = 0$ with $x = \mo O(y)$, we obtain
\begin{equation*}
	0 = \myg(t_y) + \mo O(t_y\myg(t_y)^2) - y (1 + \mo O(\myg(t_y) + t_y))
	\eqcm
\end{equation*} 
and in the setting of $x = 0$, we have
\begin{equation*}
	0 = \myg(t_y) + \mo O(t_y\myg(t_y)^2) - y (1 + \mo O(t_y\myg(t_y) + t_y^2))
	\eqfs
\end{equation*}
For $a, b, u\in \R$ with $\abs{b}\leq\frac12$, it holds
\begin{equation*}
	\abs{\frac{u+a}{1+b} - u} \leq 2\abs{a} + 2\abs{ub}\eqfs
\end{equation*}
Applied to the equations above with $u = \myg(t_y)$, $a = \mo O(t_y\myg(t_y)^2)$, and $b = \mo O(\myg(t_y) + t_y) = \mo o(1)$, this yields
\begin{equation*}
	y = \myg(t_y) + \mo O(\myg(t_y)^2 + t_y\myg(t_y))
\end{equation*}
for $x = \mo O(y)$, and for $x = 0$ with $b = \mo O(t_y \myg(t_y) + t_y^2)$, 
\begin{equation*}
	y = \myg(t_y) + \mo O(t_y\myg(t_y)^2 + t_y^2\myg(t_y))
	\eqfs
\end{equation*}
In particular, we always have
\begin{equation*}
	y = \myg(t_y) + \mo o (\myg(t_y))\eqcm
\end{equation*}
which implies
\begin{equation*}
	\myg(t_y) = y + \mo o (y)\eqfs
\end{equation*}
\subsection{\autoref{prop}}
Because of symmetry we can restrict our analysis to $y \geq 0$ and $t \geq 0$ without loss of generality.
In the proof of \autoref{lmm}, we have shown
\begin{equation*}
	y = \myg(t_y) + \mo O(\myg(t_y)^2 + t_y\myg(t_y))
\end{equation*}
for $x = \mo O(y)$, and for $x = 0$, 
\begin{equation*}
	y = \myg(t_y) + \mo O(t_y\myg(t_y)^2 + t_y^2\myg(t_y))
	\eqfs
\end{equation*}
Then, with $s := \myg(t_y)$ and $t_y = f(s)$, we have
\begin{align*}
	\frac{t_y - f(y)}{f(y)} = 
	\frac{f(s)}{f(s + \mo O(s^2+sf(s)))}-1 = 
	\mo o(1)
\end{align*}
by \textbf{(A1)} in the case of $x = \mo O(y)$, and by \textbf{(A1)'} in the case of $x = 0$,
\begin{align*}
	\frac{t_y - f(y)}{f(y)} = 
	\frac{f(s)}{f(s + \mo O(s^2f(s)+sf(s)^2))}-1 = 
	\mo o(1)
	\eqfs
\end{align*}
Hence, in both cases we get
\begin{align*}
	t_y = f(y) + \mo o (f(y))
	\eqfs
\end{align*}
Furthermore, for $t\searrow 0$,
\begin{align*}
	q(t) &= \begin{pmatrix}1 \\ t \end{pmatrix} + \mo o(t)
\end{align*}
and, thus,
\begin{equation*}
	\Pi_{\mc Q}\brOf{ \begin{pmatrix} x \\ y \end{pmatrix}} = q(t_y) = \begin{pmatrix} 1 \\ f(y) \end{pmatrix} + \mo o(f(y))
	\eqfs
\end{equation*}
\subsection{\autoref{thm}}
Note that $\argmin_{p\in\mc Q} \Ex{\normof{Z - p}^2} = \argmin_{p\in\mc Q} \normof{\Ex{Z} - p}$, as $\Ex{\normof{Z - p}^2} = \normof{\Ex{Z} - p}^2  - \normof{\Ex{Z}}^2 + \Ex{\normof{Z}^2}$.
As $\Ex{Z} = 0$, $r(0) = 1$, and $r(t) > 1$ for $t > 0$, the projected mean $m$ of $Z$ is unique and equal to $q(0)$. 

Let $(\bar X_n \rvs \bar Y_n)\pr := \bar Z_n = \frac1n \sum_{i=1}^n Z_i$. Fix $s \geq 0$. Our goal is to show \begin{equation}\label{eq:tn}
	\PrOf{t_n \leq f\brOf{\frac{s}{\sqrt{n}}}} 
	\to 
	\Phi\brOf{\frac{s}{\sigma}}\eqfs
\end{equation}
For $L, \delta > 0$ define the following events,
\begin{align*}
	A_{n,L} &:= \cb{|\bar X_n| \leq L|\bar Y_n|}\eqcm\\
	B_{n, s} &:= \cb{t_n \leq f\brOf{\frac{s}{\sqrt{n}}}}\eqcm\\
	C_{n, s} &:= \cb{\sqrt{n} \bar Y_n + \Delta_n \leq s}\eqcm\\
	D_{n, s, \delta} &:= \cb{\sqrt{n} \bar Y_n \leq s (1 + \delta)}
	\eqcm
\end{align*}
where $\Delta_n := \sqrt{n} \br{\myg(t_n) - \bar Y_n}$. Fix $\epsilon > 0$. We show  \eqref{eq:tn} by proving $\abs{\Prof{B_{n,s}} - \Phi\brOf{\frac{s}{\sigma}}} < 5\epsilon$ for $n$ large enough. We achieve this by splitting the left hand side into five parts by means of the triangle inequality and bound each summand by $\epsilon$:
\begin{enumerate}[label=(\roman*)]
	\item By the central limit theorem for $(\bar X_n, \bar Y_n)\pr$, with $\Vof Y = \sigma^2 > 0$, there is $L > 0$ and $n_1 \in \N$ such that $\Prof{A_{n, L}^{\ms c}} < \epsilon$ for all $n > n_1$.
	Thus, $\abs{\Prof{B_{n,s}} - \Prof{B_{n, s} \cap A_{n,L}}} < \epsilon$.
	\item
	Choose $\delta > 0$ such that $\abs{\Phi\brOf{\frac{s}{\sigma(1+\delta)}} - \Phi\brOf{\frac{s}{\sigma}}} + \abs{\Phi\brOf{\frac{s}{\sigma(1-\delta)}} - \Phi\brOf{\frac{s}{\sigma}}} < \epsilon$.
	\item
	By \autoref{lmm}, on the event $A_{n,L}$ for $\bar Y_n$ small enough, $\myg(t_n) = \bar Y_n + \mo o(\bar Y_n)$. Thus, there is $n_2 \in \N$ such that $\Prof{\cb{\abs{\Delta_n} > \sqrt{n} \delta \abs{\bar Y_n}} \cap A_{n, L}} \leq \epsilon$ for all $n > n_2$.
	Therefore,
	$\Prof{D_{n,s, -\delta} \cap A_{n,L}} - \epsilon < \Prof{C_{n,s} \cap A_{n,L}} < \Prof{D_{n,s, \delta} \cap A_{n,L}} + \epsilon$.
	\item
	As in (i), $\abs{\Prof{D_{n,s,\pm\delta}} - \Prof{D_{n, s,\pm\delta} \cap A_{n,L}}} < \epsilon$ for all $n > n_1$.
	\item
	By the central limit theorem, there is $n_3 \in \N$ such that $\abs{\Prof{D_{n,s,\pm\delta}} - \Phi\brOf{\frac{s}{\sigma(1\pm\delta)}}} < \epsilon$ for all $n > n_3$.
\end{enumerate}
As $B_{n, s} = C_{n,s}$, trivially  $\Prof{B_{n, s} \cap A_{n,L}} = \Prof{C_{n, s} \cap A_{n,L}}$.
All points above together yield $\abs{\Prof{B_{n,s}} - \Phi\brOf{\frac{s}{\sigma}}} < 5\epsilon$ for all $n > \max(n_1, n_2, n_3)$.
Hence, we have shown \eqref{eq:tn}.
As
\begin{equation*}
	\begin{pmatrix} m_{n,1} \\ m_{n,2} \end{pmatrix} = m_n = q(t_n) = \begin{pmatrix} 1 \\ t_n \end{pmatrix} + \mo o(t_n)\eqcm
\end{equation*}
equation \eqref{eq:tn} implies 
\begin{align*}
	\PrOf{ m_{n,2} \leq f\brOf{\frac{s}{\sqrt{n}}}} 
	&\to 
	\Phi\brOf{\frac{s}{\sigma}}\eqcm
	\\
	\PrOf{ \abs{m_{n,1}-1} \geq f\brOf{\frac{s}{\sqrt{n}}}} 
	&\to 
	0\eqfs
\end{align*}
The results for  $-t_n$ and $-m_{n,2}$  are due to symmetry.
\subsection{\autoref{coro}}
We only show the statements for $t_n$ as the results for $y_n$, $-t_n$, $-y_n$ follow similarly. Denote $F(s) :=  \Phi\brOf{\frac{s}{\sigma}}$ and let $s\geq0$.
\begin{enumerate}[label=\environmentEnumerateLabel]
\item 
	It is easy to see that \textbf{(A0)} holds for $f(y) = y^{\gamma}$. Thus, by \autoref{thm},
	\begin{equation*}
		\PrOf{n^{\frac{\gamma}{2}} t_n \leq s} = \PrOf{t_n \leq f\brOf{\frac{s^{\frac1\gamma}}{\sqrt{n}}}} \xrightarrow{n\to\infty} F\brOf{s^\frac1\gamma}
		\eqfs
	\end{equation*}
	Furthermore, $r(t) = 1 + \int_0^t x^{\frac1\gamma} \dl x = 1+\frac{\gamma}{1 + \gamma} t^{\frac{1 + \gamma}{\gamma}}$.
\item
	It is easy to check \textbf{(A0)} for $f(y) = \br{-\log(y)}^{-\gamma}$.
	
	The inverse function of $f$ is $\myg(x) = \exp\brOf{-x^{-\frac1\gamma}}$, which yields the expression for $r(t)$. By \autoref{thm},
	\begin{equation*}
		\PrOf{t_n \leq f\brOf{\frac{s}{\sqrt{n}}}} \xrightarrow{n\to\infty} F(s)\eqfs
	\end{equation*}
	It holds
	\begin{equation*}
		\PrOf{t_n \leq f\brOf{\frac{s}{\sqrt{n}}}} =
		\PrOf{\br{\frac12\log(n)}^{\gamma} t_n \leq \br{\frac{\log(\sqrt{n})}{\log\brOf{\frac{\sqrt{n}}{s}}}}^\gamma}
		\eqfs
	\end{equation*}
	As $\log(\sqrt{n})/\log(\sqrt{n}/s) \xrightarrow{n\to\infty} 1$ for all $s > 0$, 
	\begin{align*}
		\PrOf{\br{\frac12\log(n)}^{\gamma} t_n \leq t} \xrightarrow{n\to\infty}
		\begin{cases}
		F(0) = \frac12\quad\text{for } 0 < t < 1\eqcm\\
		F(\infty) = 1\quad\text{for } t > 1\eqcm\\
		\end{cases}
	\end{align*}
	which, together with symmetry of the distribution, shows convergence of $\br{\frac12\log(n)}^{\gamma} t_n$ in distribution to a uniform distribution on $\cb{-1, 1}$.
\item
	It is easy to check \textbf{(A0)} for $f(y) = \exp\brOf{-y^{-\gamma}}$.	
	The inverse function of $f$ is $\myg(x) = (- \log(x))^{-\frac1\gamma}$, which yields the expression for $r(t)$.
	
	Let $c, u > 0$. For $s\in(1,\infty)$ and $n$ large enough, $u \exp(-(\sqrt{n}/c)^{\gamma}) \leq \exp(-(\sqrt{n}/(cs))^\gamma) = f(cs n^{-\frac12})$. Thus,
	with $U_{n,c} := \exp\brOf{(\sqrt{n}/c)^\gamma} t_n$, 
	\begin{equation*}
		\PrOf{U_{n,c} \leq u} \leq \PrOf{t_n \leq f\brOf{\frac{cs}{\sqrt{n}}}} \xrightarrow{n\to\infty}F(cs)
	\end{equation*}
	by \autoref{thm}.
	Similarly, for $s \in (0,1)$, 
	\begin{equation*}
		\PrOf{U_{n,c} \leq u} \geq \PrOf{t_n \leq f\brOf{\frac{cs}{\sqrt{n}}}} \xrightarrow{n\to\infty} F(cs)
		\eqfs
	\end{equation*}
	Thus, 
	\begin{equation*}
		\PrOf{U_{n,c} \leq u} \xrightarrow{n\to\infty} F(c) =: p_c\eqcm
	\end{equation*}
	which implies $\Prof{U_{n,c} \geq u} \xrightarrow{n\to\infty} 1-p_c$. As $t_n$ is symmetric, $\Prof{U_{n,c} \leq -u} \xrightarrow{n\to\infty} 1-p_c$, which leaves $\PrOf{\abs{U_{n,c}} < u} \xrightarrow{n\to\infty} 2p_c-1$.
\end{enumerate}
\subsection{\autoref{rem2}}\label{ssec:proof:rem2}
\begin{enumerate}[label=\environmentEnumerateLabel]
\item 
	It is easy to see that \textbf{(A1)} hold for $f(y) = y^{\gamma}$.
\item
	To verify \textbf{(A1)} for $f(y) = \br{-\log(y)}^{-\gamma}$, note
	\begin{equation*}
		\lim_{y\to 0}\frac{\log(y)}{\log(y+h(y))} = \lim_{y\to 0} \frac{y + y h\pr(y)}{y + h(y)} = 1
	\end{equation*}
	for $h(y) = \mo o(y)$. Here we use $h(y) = c y\br{y + \log(\frac1y)^{-\gamma}}$.
\item
	To verify \textbf{(A1)'} for $f(y) = \exp\brOf{-y^{-\gamma}}$, note
	\begin{equation*}
		\frac{\exp\brOf{-(y+h(y))^{-a}}}{\exp\brOf{-y^{-a}}} 
		=
		\exp\brOf{y^{-a}-(y+h(y))^{-a}}
		\xrightarrow{y\to0} 1
	\end{equation*}
	for $a >0$ and $h(y) = \mo o(y^2)$, as
	\begin{equation*}
		y^{-a}-(y+h(y))^{-a} \to 0
		\eqfs
	\end{equation*}
	Here, we use $h(y) = c \br{y^2 \exp\brOf{-y^{-\gamma}} +  y \exp\brOf{-2y^{-\gamma}}}$.
\end{enumerate}
\section*{Acknowledgments}
The author gratefully  acknowledges  support  by  the  German  Research  Foundation (DFG) through the Research Training Group RTG 1953.
Furthermore, the author is thankful to Benjamin Eltzner and Stephan Huckeman for motivation to write this paper, a fruitful discussion of the topic, and many helpful comments on drafts of this paper. Additionally, the author thanks Sandra Schluttenhofer and Jan Johannes for reading and commenting on drafts of the paper.
\bibliographystyle{alpha}
\bibliography{literature}
\end{document}